\documentclass[12pt]{article}
\usepackage{prelude}

\newkeytheorem{theorem}[name = Theorem, style = theorem]
\newkeytheorem{lemma}[name = Lemma, style = theorem]
\newkeytheorem{proof}[name = Proof, style = proof, numbered=no]

\newcommand{\Bu}{\textrm B}

\title{A general identity for umbral operators and a special subclass}
\author{Kei Beauduin}
\date{}

\begin{document}

\maketitle

\begin{abstract}
    We prove a new universal identity for umbral operators. This motivates the definition of a subclass satisfying a simplified identity, which we fully characterize. The results are illustrated with common examples of the theory of umbral calculus.
\end{abstract}

%\noindent\textbf{AMS Classifications :} 05A40, 13F25, 11B73, 05A10, 05A19, 47B99, 11B68

\section{Introduction}

In \cite{beauduin2025}, we derived an operational identity that is instrumental for proving the general umbral identity of \Cref{t:phixn}. This identity, along with numerous examples, motivated the definition of a subclass of umbral operators satisfying a simpler relation, which we characterize fully in \Cref{t:special}. Several such examples are detailed in \Cref{s:ex}.

\section{Proof of the identity}

We recall some umbral calculus notation \cite{beauduin2024,rota1973}. An \emph{umbral operator} $\phi$ is the linear operator mapping monomials to the polynomial sequence $\{\phi_n(x)\}_{n\in\N}$ determined by the generating function
\begin{equation}\label{e:gen}
    e^{xf(t)} = \sum_{n=0}^\infty \phi_n(x) \frac{t^n}{n!},
\end{equation}
for some power series $f$ satisfying $f(0) = 0$ and $f'(0) \neq 0$, i.e., $f$ admits a compositional inverse $f^{-1}$.
We define the \emph{coefficients} $\coeff{n}{k}_\phi$ of $\phi$ by
\begin{equation}\label{e:coeff}
    \phi x^n = \phi_n(x) = \sum_{k=0}^n \coeff{n}{k}_\phi x^k.
\end{equation}
The \emph{partial Bell polynomials} $B_{n,k}$ \cite[Chap.~3.3]{comtet1974} are defined by
\begin{equation}\label{e:Bell}
    \frac{f(t)^k}{k!} = \sum_{n=k}^\infty B_{n,k}(f^{(1)}(0), \dots, f^{(n-k+1)}(0)) \frac{t^n}{n!},
\end{equation}
where $f^{(j)}$ denotes the $j$-th derivative of $f$. Hence \cref{e:gen,e:coeff} imply
\begin{equation}\label{e:coeffBell}
    \coeff{n}{k}_\phi = B_{n,k}(f^{(1)}(0), \dots, f^{(n-k+1)}(0)).
\end{equation}
We define the \emph{delta operator} associated to $\phi$ by
\begin{equation}\label{e:Q}
    Q := f^{-1}(D),
\end{equation}
where $D$ is the operator of differentiation with respect to $x$. We have
\begin{equation}\label{e:Qphi}
    Q \phi_n(x) = n \phi_{n-1}(x)
\end{equation}
for all positive integers $n$. Let $\x$ denote multiplication by $x$. We can now state the main identity.
\begin{theorem}\label{t:phixn}
    For all nonnegative integers $n$,
    \begin{equation}\label{e:phixn}
        \phi \x^n = \sum_{k=0}^n \x^k  B_{n, k}(f^{(1)}(Q), \ldots, f^{(n-k+1)}(Q)) \phi.
    \end{equation}
\end{theorem}
When applied to $1$, both sides of \cref{e:phixn} reduce to $\phi_n(x)$ by \cref{e:coeff,e:coeffBell}. At $n=1$, \Cref{t:phixn} reduces to the identity
\begin{equation}\label{e:recu}
    \phi \x = \x (Q')^{-1} \phi,
\end{equation}
which is the operational analog of the recurrence relationship $\phi_{n+1}(x) = x (Q')^{-1} \phi_n(x)$ \cites[Thm.~4 (4)]{rota1973}{beauduin2024}. \Cref{e:phixn} can also be written as
\begin{equation}\label{e:phixn2}
    \phi \x^n = \sum_{k=0}^n \x^k \phi B_{n, k}(f^{(1)}(D), \ldots, f^{(n-k+1)}(D))
\end{equation}
as a consequence of $Q \phi = \phi D$, which itself is the operational reformulation of \cref{e:Qphi}. We prove \Cref{t:phixn} using the following lemma.

\begin{lemma}[\cite[Cor.~4.1]{beauduin2025}]\label{l:V}
    For a power series $V$, if $h$ is a function such that $h'(t) = 1/V(t)$, then for all $s\in\C$
    \[
    e^{\x V(D)s} = (e^\x)^{h^{-1}(h(D) + s) - D}.
    \]
\end{lemma}
The exponentiation of operators used in \Cref{l:V} was introduced in \cite{beauduin2025a} by
\[
A^B = \sum_{k=0}^\infty (A-1)^k \frac{(B)_k}{k!},
\]
for operators $A, B$ and where $(B)_k := B(B-1)\dots(B-k+1)$ is the \emph{falling factorial}. This exponentiation satisfies the identity
\begin{equation}\label{e:eAB}
    (e^A)^B = \sum_{k=0}^\infty \frac{A^k B^k}{k!}.
\end{equation}

\begin{proof}[of \Cref{t:phixn}]
    It follows from \cref{e:recu} that
    \begin{equation}\label{e:nrecu}
        \phi \x^n =\x (Q')^{-1} \phi \x^{n-1} = \ldots = (\x (Q')^{-1})^n \phi,
    \end{equation}
    for every $n\in\N$. Multiplying \cref{e:nrecu} by $t^n/n!$ and summing over $n$ yields the equivalent generating function identity
    \begin{equation}\label{e:phiexp}
        \phi e^{\x t} = e^{\x (Q')^{-1} t} \phi,
    \end{equation}
    and according to \Cref{l:V},
    \begin{equation}\label{e:expQ'-1}
        e^{\x (Q')^{-1}t} = (e^\x)^{f(Q + t) - D}.
    \end{equation}
    By the Taylor expansion
    \[
    f(Q + t) - D = \sum_{n=0}^\infty f^{(n)}(Q) \frac{t^n}{n!} - D =  \sum_{n=1}^\infty f^{(n)}(Q) \frac{t^n}{n!},
    \]
    thus, by \cref{e:eAB} and (\ref{e:Bell})
    \[
    (e^\x)^{f(Q + t) - D} = \sum_{k=0}^\infty \frac{\x^k}{k!} \pa{\sum_{n=1}^\infty f^{(n)}(Q) \frac{t^n}{n!}}^k = \sum_{k=0}^\infty \x^k \sum_{n=k}^\infty B_{n, k}(f^{(1)}(Q), \ldots, f^{(n-k+1)}(Q))  \frac{t^n}{n!},
    \]
    and interchanging the order of summation yields
    \begin{equation}\label{e:Bellproof}
        (e^\x)^{f(Q + t) - D} = \sum_{n=0}^\infty \frac{t^n}{n!} \sum_{k=0}^n \x^k  B_{n, k}(f^{(1)}(Q), \ldots, f^{(n-k+1)}(Q)).
    \end{equation}
    We combine \cref{e:phiexp,e:expQ'-1,e:Bellproof} to obtain
    \[
    \phi e^{\x t} = \sum_{n=0}^\infty \frac{t^n}{n!} \sum_{k=0}^n \x^k  B_{n, k}(f^{(1)}(Q), \ldots, f^{(n-k+1)}(Q)) \phi,
    \]
    and by identifying the coefficient in front of $t^n/n!$, we arrive at \cref{e:phixn}.
\end{proof}

\section{A special class of umbral operators}

In this section, we characterize the umbral operators satisfying the identity
\begin{equation}\label{e:special}
    \phi \x^n = \sum_{k=0}^n \coeff{n}{k}_\phi \x^k U^k V^n \phi,
\end{equation}
where $U$ and $V$ are power series in $D$. By uniqueness of the normal-ordered form \cite{mansour2015}, \Cref{t:phixn} then shows that \cref{e:special} is equivalent to
\[
B_{n, k}(f^{(1)}(Q), \ldots, f^{(n-k+1)}(Q)) = \coeff{n}{k}_\phi U^k V^n.
\]

The first step is to characterize the functions satisfying a generalized form of Cauchy's classical functional equation.

\begin{lemma}\label{l:funcEq}
    Let $g$ be an invertible formal power series of the form
    \[
    g(x) = \sum_{k=1}^\infty \frac{c_k}{k!} x^k, \quad c_1 \neq 0.
    \]
    There exist formal power series $u, v$ such that $g$ satisfies the formal power series functional equation
    \begin{equation}\label{e:funcEq}
        g(x+y) = g(x) + u(x) g(v(x) y)
    \end{equation}
    if and only if $g$ has one of the forms listed in \cref{e:table1} below. The corresponding choices of $u$ and $v$ are also listed.
    \begin{equation}\label{e:table1}
        \begin{array}{c|c|c}
        g(x) & u(x) & v(x) \\
        \hdoubleline
        ax & A(x) & A(x)^{-1} \\
        \hline
        \rule{0pt}{10pt}
        a(e^{bx} - 1) & e^{bx} & 1 \\
        \hline
         a\log(1+bx) & 1 & (1+bx)^{-1} \\
         \hline
         a((1+bx)^c - 1) & (1+bx)^c & (1+bx)^{-1}
        \end{array}
    \end{equation}
    provided $a, b, c \neq 0$, and where $A$ is a power series satisfying $A(0) \neq 0$.
\end{lemma}

\begin{proof}
    After differentiating \cref{e:funcEq} once and twice at $y = 0$, we find that
    \[
    g'(x) = u(x)v(x)c_1, \qquad g''(x) = u(x)v(x)^2c_2.
    \]
    If $c_2 = 0$, then $g''(x) = 0$ and with the additional fact that $g(0) = 0$, we infer that $g(x) = c_1 x$. For this function, the condition $u(x)v(x) = 1$ is necessary to satisfy \cref{e:funcEq}. If $c_2 \neq 0$, we can solve for $u$ and $v$:
    \[
    u(x) = \pa{\dfrac{g'(x)}{c_1}}^2 \dfrac{c_2}{g''(x)}, \qquad v(x) = \dfrac{g''(x)}{c_2} \dfrac{c_1}{g'(x)}.
    \]
    \Cref{e:funcEq} can be rearranged, and by replacing $y$ by $y g'(x)/c_1$, we get
    \[
    \frac{g''(x)}{c_2} \pa{g\pa{x + \frac{g'(x)}{c_1} y} - g(x)} = \pa{\frac{g'(x)}{c_1}}^2 g\pa{\frac{g''(x)}{c_2} y}.
    \]
    We now differentiate at $x = 0$:
    \[
    \frac{c_3}{c_2} g(y) + \pa{1+ \frac{c_2}{c_1}y} g'(y) - c_1 = 2 \frac{c_2}{c_1} g(y) + \frac{c_3}{c_2}y g'(y).
    \]
    Therefore, we only need to solve the linear differential equation in $y$
    \[
    (1+ \alpha y)g'(y) = \beta g(y) + c_1, \quad g(0) = 0,
    \]
    where
    \[
    \alpha := \frac{c_2}{c_1} - \frac{c_3}{c_2},\qquad \beta := 2 \frac{c_2}{c_1} - \frac{c_3}{c_2}.
    \]
    There are four cases, but the (linear) case $\alpha = \beta = 0$ was already covered.
    \begin{itemize}
        \item If $\alpha = 0$ and $\beta \neq 0$, then $g(y) = c_1 \dfrac{e^{\beta y}-1}{\beta}$.
        \item If $\alpha \neq 0$ and $\beta = 0$, then $g(y) = c_1 \dfrac{\log(1+\alpha y)}{\alpha}$.
        \item If $\alpha \neq 0$ and $\beta \neq 0$, then $g(y) = \dfrac{c_1}{\beta}((1+\alpha y)^{\beta/\alpha}-1)$.
    \end{itemize}
    These candidate functions satisfy \cref{e:funcEq} and are precisely the forms listed in \cref{e:table1}, which concludes the proof.
\end{proof}

Note that the set of solutions is closed under inversion, that is, if $g$ is a solution of \cref{e:funcEq}, then its compositional inverse $g^{-1}$ also satisfies an equation of the form \cref{e:funcEq}, for suitable choices of $u$ and $v$.

For $a\in\C$, let $E^a$ be the shift operator defined by $E^a p(x) = p(x+a)$ for a polynomial $p$. We recall that an operator $U$ is said to be shift-invariant whenever it commutes with $E^a$ for all $a$. Shift-invariance is equivalent to the existence of a formal power series $\tilde U$ such that $\tilde U(D) = U$ \cite{rota1973,beauduin2024}, so by definition (\ref{e:Q}), $Q$ is an example of such an operator.

\begin{theorem}\label{t:special}
    Let $\phi$ be an umbral operator with the associated delta operator $Q$. $\phi$ satisfies, for some shift-invariant operators $U, V$, the commutation identity \textup{(\ref{e:special})} if and only if $Q$ has one of the forms listed in \cref{e:table2} below. The corresponding choices of $U$ and $V$ are also listed.
    \begin{equation}\label{e:table2}
        \begin{array}{c|c|c}
        Q & U & V \\
        \hdoubleline
        aD & A & A^{-1} \\
        \hline
        \rule{0pt}{10pt}
        a (e^{bD}-1) & 1 & e^{-bD} \\
        \hline
        a \log(1+bD) & 1+bD & 1 \\
        \hline
        a((1+bD)^c-1) & 1+bD & (1+bD)^{-c}
        \end{array}
    \end{equation}
    provided $a, b, c \neq 0$, and where $A$ is a shift-invariant invertible operator.
\end{theorem}

When $\phi$ is associated to the delta operator $Q = f^{-1}(D)$, $\phi^{-1}$ is associated to $f(D)$ \cite{beauduin2024}. It follows directly from \Cref{t:special} that this class is stable under inversion.

\begin{proof}
    According to \cref{e:phiexp,e:expQ'-1}
    \begin{equation}\label{e:phiexp2}
        \phi e^{\x t} = (e^\x)^{f(Q+t)-D} \phi.
    \end{equation}
    Since $U$ and $V$ are shift-invariant, they can be expressed as formal power series in $D$: $\tilde U(D) = U$ and $\tilde V(D) = V$. Defining $u := \tilde U \circ f$ and $v := \tilde V \circ f$, the operational equation
    \begin{equation}\label{e:opeq}
        f(Q + t) = D + U f(Vt)
    \end{equation}
    is equivalent, via the isomorphism theorem \cite[Thm.~3]{rota1973}, to the functional equation
    \[
    f(x + t) = f(x) + u(x) f(v(x)t).
    \]
    Thus, by \Cref{l:funcEq}, for $a, b, c \neq 0$, $f(t)$ is either $at$, $a(e^{bt}-1)$, $a\log(1+bt)$ or $a((1+bt)^c - 1)$. Consequently, the delta operator $Q = f^{-1}(D)$ takes the form listed in \cref{e:table2}. The corresponding expressions for $U = u(Q)$ and $V = v(Q)$ are then computed directly. By \cref{e:Bell,e:coeffBell}
    \begin{align}\label{e:exp=UV}
        &(e^\x)^{Uf(Vt)} = \sum_{k=0}^\infty \frac{\x^k}{k!} (U f(Vt))^k \nonumber\\
        ={}& \sum_{k=0}^\infty \x^k U^k \sum_{n=k}^\infty \coeff{n}{k}_\phi \frac{(Vt)^n}{n!} = \sum_{n=0}^\infty \frac{t^n}{n!} \sum_{k=0}^n \coeff{n}{k}_\phi \x^k U^k V^n.
    \end{align}
    Combining \cref{e:phiexp2,e:opeq,e:exp=UV}, we compute
    \[
    \phi e^{\x t} = (e^\x)^{f(Q+t)-D} \phi = (e^\x)^{Uf(Vt)} \phi = \sum_{n=0}^\infty \frac{t^n}{n!} \sum_{k=0}^n \coeff{n}{k}_\phi \x^k U^k V^n \phi.
    \]
    We identify the coefficients in front of $t^n/n!$ to recover the desired identity. Conversely, summing \cref{e:special} over $n$ gives the same generating-function identity, and the uniqueness of the normal-ordered form recovers \cref{e:opeq}; hence the preceding steps are reversible.
\end{proof}

\section{Some examples}\label{s:ex}

First note that $E^a=e^{aD}$, the operator form of Taylor expansion. When the delta operator $Q$ is the \emph{forward difference} operator $\Delta:= E^1 - 1$, we denote the associated umbral operator by $\varphi$ that maps $x^n$ to the falling factorial $(x)_n$. According to \Cref{t:special}, in that case
\begin{equation}
    \varphi \x^n = \sum_{k=0}^n \coeff{n}{k}_\varphi \x^k E^{-n} \varphi = (\x)_n E^{-n} \varphi,
\end{equation}
where $\coeff{n}{k}_\varphi$ is the signed \emph{Stirling number of the first kind} \cite{jordan1933}. If we apply this operational equality to a monomial $x^m$, we obtain
\begin{equation}
    (x)_{n+m} = (x)_n (x-n)_m,
\end{equation}
which follows directly from the definition of falling factorials.

\medskip

The delta operator associated to $\varphi^{-1}$ is $\log(1+D)$, and therefore, \Cref{t:special} gives
\begin{equation}\label{e:Touchard}
    \varphi^{-1} \x^n = \sum_{k=0}^n \coeff{n}{k}_{\varphi^{-1}} \x^k (1+D)^k \varphi^{-1} = \sum_{k=0}^n \coeff{n}{k}_{\varphi^{-1}} \x^k \varphi^{-1} E^k,
\end{equation}
where the second formula is in the form of \cref{e:phixn2}. Because $\varphi^{-1} x^n$ corresponds to the Touchard polynomial $T_n(x)$ \cite{touchard1956,beauduin2024} and $\coeff{n}{k}_{\varphi^{-1}}$ to the \emph{Stirling number of the second kind} $\Stir{n}{k}$ \cite{jordan1933}, applying \cref{e:Touchard} to $x^m$ yields
\begin{equation}
    T_{n+m}(x)
    = \sum_{k=0}^n \Stir{n}{k} \x^k \varphi^{-1} (x+k)^m
    = \sum_{k=0}^n \Stir{n}{k} x^k \sum_{j=0}^m \binom{m}{j} k^{m-j} T_j(x),
\end{equation}
which is a known extension of Spivey's identity for \emph{Bell numbers} $B_n := T_n(1)$ \cite{gould2008,spivey2008,beauduin2024}.

\medskip

\Cref{t:special} also applies to the \emph{Laguerre umbral operator} $L$, associated with $D/(1-D) = (1-D)^{-1} - 1$, which maps the $n$-th monomial to the \emph{Laguerre polynomial} $L_n(x)$:
\begin{equation}
    L \x^n = \sum_{k=0}^n \coeff{n}{k}_L \x^k (1-D)^{k+n} L,
\end{equation}
where $\coeff{n}{k}_L$ are the signed \emph{Lah numbers} $\binom{n-1}{k-1}\frac{n!}{k!} (-1)^{n-k}$ \cite{lah1954,lah1955}. Applying this to $x^m$ yields the polynomial identity
\begin{equation}
    L_{n+m}(x) = \sum_{k=0}^n \coeff{n}{k}_L x^k L^{(k+n)}_m(x),
\end{equation}
where $L^{(\alpha)}_m(x) := (1-D)^\alpha L_m(x)$ is an \emph{associated Laguerre polynomial} \cite{beauduin2024,comtet1974,rota1973,dibucchianico1997}.

\medskip

Now consider the delta operator $D(1-D)$ associated to the \emph{Bucchianico polynomials} $\Bu_n(x)$ \cite{beauduin2024}, introduced by Di Bucchianico \cite[Ex.~3.2.7 c]{dibucchianico1997}. Since
\[
D(1-D) = \frac{1 - (1-2D)^2}{4},
\]
\Cref{t:special} applies and yields
\begin{equation}
    \Bu \x^n = \sum_{k=0}^n \coeff{n}{k}_\Bu \x^k  (1-2D)^{k-2n} \Bu,
\end{equation}
where
\begin{equation}
    \coeff{n}{k}_\Bu = \binom{2n-k-1}{n-1} \frac{(n-1)!}{(k-1)!},
\end{equation}
for $1 \le k \le n$ and equal to $1$ if $n=k=0$ and $0$ otherwise (see \cite{dibucchianico1997,beauduin2024}).

\printbibliography

\end{document}